\newif\ifpictures
\picturestrue 

\documentclass[11pt]{amsart}
\usepackage{amssymb}
\usepackage{epsf}

\numberwithin{equation}{section}

\newtheorem{thm}{Theorem}
\newtheorem{prop}[thm]{Proposition}

\newtheorem{example}[thm]{Example}
\newtheorem{remark}[thm]{Remark}

\newenvironment{rem}{\begin{remark}\rm}{\end{remark}}

\newcounter{FNC}[page]
\def\silentfootnote#1{{\addtocounter{FNC}{2}$^\fnsymbol{FNC}$%
     \let\thefootnote\relax\footnotetext{$^\fnsymbol{FNC}$#1}}}

\newcommand{\hh}{\hspace{8.525pt}}  

\newcommand{\bc}{{c}}
\newcommand{\be}{{e}}
\newcommand{\bp}{{p}}
\newcommand{\bv}{{v}}

\title{Lines Tangent to $2n-2$ spheres in ${\mathbb R}^n$}

\author{Frank Sottile}
\address{Department of Mathematics\\
        University of Massachusetts\\
        Lederle Graduate Research Tower\\
        Amherst, Massachusetts, 01003\\
        USA}
\email{sottile@math.umass.edu}
\urladdr{http://www.math.umass.edu/\~{}sottile}

\author{Thorsten Theobald}
\address{Zentrum Mathematik\\
        Technische Universit\"at M\"unchen\\
        D--80290 M\"unchen\\
        Germany}
\email{theobald@mathematik.tu-muenchen.de}
\urladdr{http://www-m9.mathematik.tu-muenchen.de/\~{}theobald/}

\thanks{Research of first author supported in part by NSF grant DMS-0070494}

\copyrightinfo{}{}

\begin{document}

\begin{abstract}
We show that for $n \ge 3$
there are $3 \cdot 2^{n-1}$ complex com\-mon tangent lines to
$2n-2$ general spheres in $\mathbb{R}^n$ and that there is a choice
of spheres with all common tangents real.
\end{abstract}

\maketitle

\section{Introduction}

We study the following problem from (real) enumerative geometry.

\medskip

\begin{description}
\item[Given] $2n-2$ (not necessarily disjoint)
  spheres with centers $c_i \in \mathbb{R}^n$ 
  and radii $r_i$, $1 \le i \le 2n-2$. \smallskip

\item[Question] In the case of finitely many common tangent lines,
what is their maximum number? 
\end{description}

The number of $2n-2$ spheres guarantees that in the generic case there is
indeed a finite number of common tangent lines. In
particular, for $n=2$ the answer is~4 since two disjoint circles have 4 common
tangents. 

The reason for studying this question---which, of course, is an appealing 
and fundamental geometric question in itself---came from different 
motivations.
An essential task in statistical analysis is to find the line that
best fits the data in the sense of minimizing the maximal distance to
the points~(see, e.g., \cite{chan-2000}).
More precisely, the decision variant of this problem asks:
Given $m,n \in \mathbb{N}$, $r > 0$, and a set of points
$y_1, \ldots, y_m \in \mathbb{R}^n$,
does there exist a line $l$ in $\mathbb{R}^n$ such every point $y_i$ has
Euclidean
distance at most $r$ from~$l$.
From the complexity-theoretical point of view,
for fixed dimension the problem can be solved in polynomial time via
quantifier elimination over the reals~\cite{fks-96}.
However, currently no practical algorithms focusing on exact computation
are known for $n > 3$ (for approximation algorithms,
see~\cite{chan-2000}).

From the algebraic perspective, for dimension~3 it was shown 
in~\cite{aas-99,ssty-2000} how to reduce the algorithmic problem to an 
algebraic-geometric core problem: finding
the real lines which all have the same prescribed distance from 4~given points;
or, equivalently, finding the real common tangent lines to 4~given unit 
spheres in $\mathbb{R}^3$. 
This problem in dimension 3  was treated in~\cite{MPT01}.

\begin{prop}\label{prop:MPTh}
Four unit spheres in $\mathbb{R}^3$ have at most $12$ common tangent lines 
unless their centers are collinear.
Furthermore, there exists a configuration with $12$ different real tangent
lines.
\end{prop}

The same reduction idea to the algebraic-geometric core problem also 
applies to arbitrary dimensions, in this case leading to the 
general problem stated at the beginning.

{}From the purely algebraic-geometric point of view, this tangent 
problem is interesting for
the following reason. In dimension~3, the formulation of the problem
in terms of Pl\"ucker coordinates gives 5 quadratic equations in
projective space $\mathbb{P}^5_{\mathbb R}$, whose common zeroes in
$\mathbb{P}^5_{\mathbb C}$ include a 1-dimensional
component at infinity (accounting for the ``missing'' $2^5-12 = 20$ solutions).
Quite remarkably, as observed in \cite{aluffi-fulton-2001},
this excess component cannot be resolved by a single blow-up.
Experimental results in~\cite{sottile-macaulay-2001} for $n=4,5$, and $6$,
indicate that for higher dimensions the generic number of solutions differs
from the B\'{e}zout number of the straightforward polynomial formulation
even more. 
We discuss this further in Section~5.

Our main result can be stated as follows.


\begin{thm}   \label{th:ndimnumber}
 Suppose $n\geq 3$.
 \begin{enumerate}
 \item[(a)] Let $c_1,\ldots,c_{2n-2}\in\mathbb{R}^n$ 
            affinely span ${\mathbb R}^n$, and let $r_1,\ldots,r_{2n-2}>0$. 
            If the $2n-2$ spheres with centers $c_i$ and radii $r_i$ 
            have only a finite number of complex common tangent lines, 
            then that number is bounded by $3 \cdot 2^{n-1}$.

 \item[(b)] There exists a configuration with $3 \cdot 2^{n-1}$
            different real common tangent lines.
            Moreover, this configuration can be achieved with unit spheres.
 \end{enumerate}
\end{thm}

Thus the bound for real common tangents equals the ({\it a priori greater})
bound for complex common tangents; so this problem of common tangents to
spheres is fully real in the sense of enumerative real algebraic
geometry~\cite{So97c,So-DIMACS}.
We prove Statement (a) in Section 2 and Statement (b) in Section 3, where 
we explicitly describe configurations with $3 \cdot 2^{n-1}$
common real tangents.
Figure~\ref{Fig1} shows a configuration of 4 spheres in ${\mathbb R}^3$ with 12
common tangents (as given in~\cite{MPT01}).
\ifpictures
\begin{figure}[htb]
$$
  \epsfxsize=3.025in \epsfbox{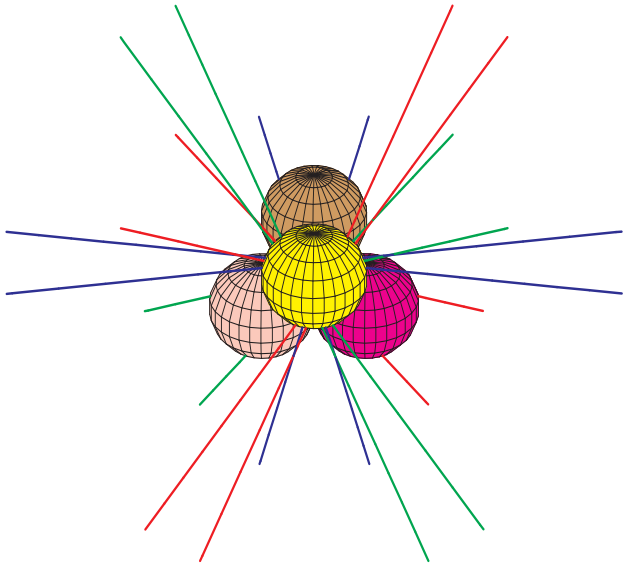}
$$
\caption{Spheres with 12 real common tangents}\label{Fig1} 
\end{figure}
\fi

In Section~4, we show that there are configurations of spheres with
affinely dependent centers having $3 \cdot 2^{n-1}$ \emph{complex}
common tangents; thus, the upper bound of 
Theorem~\ref{th:ndimnumber} also holds for spheres in this special position.
Megyesi~\cite{Me02} has recently shown that all  $3 \cdot 2^{n-1}$
may be real.
We also show that if the centers of the spheres are the vertices of the
crosspolytope in ${\mathbb R}^{n-1}$, there will be at most $2^n$ common 
tangents, and if the spheres overlap but do not contain the centroid of the
crosspolytope, then all $2^n$ common tangents will be real.
We conjecture that when the centers are affinely dependent and all spheres
have the same radius, then there will be at most $2^n$ real common tangents.
Strong evidence for this conjecture is provided by Megyesi~\cite{Me01}, who
showed that there are at most 8 real common tangents to 4 unit spheres 
in ${\mathbb R}^3$ whose centers are coplanar but not collinear.

In Section 5, we put the tangent problem into the perspective of common
tangents to general quadric hypersurfaces. In particular, we discuss 
the problem of common tangents to $2n-2$ smooth quadrics in 
\emph{projective} $n$-space, and describe the excess
component at infinity for this problem of spheres.
In this setting, Theorem~\ref{th:ndimnumber}(a) implies that there will be at
most $3 \cdot 2^{n-1}$ isolated common tangents to $2n-2$ quadrics in 
projective $n$-space, when the quadrics all contain the same (smooth) quadric
in a given hyperplane. In particular, the problem of the spheres can be
seen as the case when the common quadric is at infinity and contains no 
real points.
We conclude with the question of how many of these common tangents may be
real when the shared quadric has real points.
For $n=3$, there are 5 cases to consider, and for each, all 12 lines can be
real~\cite{sottile-macaulay-2001}. 
Megyesi~\cite{Me02} has recently shown that all common
tangents may be real, for many cases of the shared quadric.

\section{Polynomial Formulation with Affinely Independent Centers}\label{sec:indep}

For $x, y \in \mathbb{C}^n$, let $x \cdot y := \sum_{i=1}^n x_i y_i$
denote their Euclidean dot product. We write $x^2$ for $x \cdot x$.

We represent a line in ${\mathbb C}^n$ by a point $\bp\in{\mathbb C}^n$
lying on the line and a direction vector 
$\bv\in{\mathbb P}_\mathbb{C}^{n-1}$ of that line. (For notational convenience
we typically work with a representative of the direction vector in 
$\mathbb{C}^n \setminus \{0\}$.) 
If $v^2 \neq 0$ we can make $p$ unique by requiring that
$p \cdot v = 0$.

By definition, a line $\ell = (\bp,\bv)$ is tangent to the sphere with 
center $\bc\in{\mathbb R}^n$ and radius $r$
if and only if it is tangent to the quadratic
hypersurface $(x-c)^2 = r^2$, i.e., 
if and only if the quadratic equation $(p +tv - c)^2 = r^2$ has a solution of
multiplicity  two. 
When $\ell$ is real then this is equivalent to the metric
property that $\ell$ has Euclidean distance $r$ from $c$.
\ifpictures
$$
  \setlength{\unitlength}{.9167pt} 
  \begin{picture}(106,56)
   \put(0,0){\epsfxsize=91.67pt\epsfbox{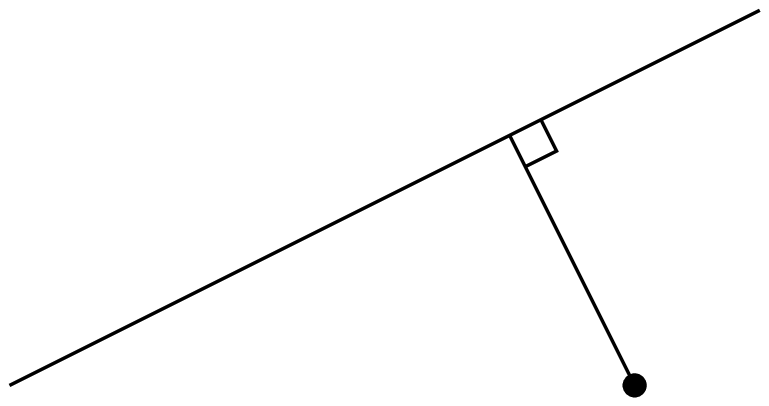}}  
   \put(28,24){$\ell$}
   \put(87,3){$\bc$}
   \put(65,16){$r$}
 \end{picture}
$$
\fi
For any line $\ell \subset \mathbb{C}^n$, the algebraic tangent condition on
$\ell$ gives the equation
$$
   \frac{[\bv\cdot(\bp-\bc)]^2}{\bv^2} - (\bp-\bc)^2\ + r^2\ =\ 0 \,.
$$
For $v^2 \neq 0$ this is equivalent to
 \begin{equation}\label{eq:tanSphere}
  \bv^2\bp^2-2\bv^2 \bp\cdot\bc + \bv^2\bc^2 - [\bv\cdot\bc]^2-r^2\bv^2\ =\ 0\,.
 \end{equation}

To prove part (a) of Theorem~\ref{th:ndimnumber}, we can choose
$\bc_{2n-2}$ to be the origin and set $r:=r_{2n-2}$. 
Then the remaining centers span ${\mathbb R}^n$.
Subtracting the equation for the sphere centered at the origin from the
equations for 
the spheres $1,\ldots,2n-3$ gives the system
 \begin{equation}\label{eq:origintransform}
  \begin{array}{rcl}
    \bp\cdot\bv &=& 0 \, ,\\ \rule{0pt}{20pt}
    \bp^2 &=& r^2 \, , \quad\mbox{and} \\ \rule{0pt}{20pt}
      2 \bv^2\bp\cdot\bc_i & = & \bv^2 \bc_i^2 - [\bv \cdot \bc_i]^2 
      - \bv^2 (r_i^2 - r^2) \, ,
      \qquad i=1,2,\ldots,2n{-}3\,.
  \end{array}
 \end{equation}
\begin{rem}
 Note that this system of equations does not have a solution with
 $v^2 = 0$. Namely, if we had $v^2 = 0$, then 
 $\bv \cdot \bc_i=0$ for all $i$.
 Since the centers span ${\mathbb R}^n$, this would imply $\bv=0$, 
 contradicting $\bv\in{\mathbb P}^{n-1}_{\mathbb C}$.
 This validates our assumption that $v^2\neq 0$ prior to~(\ref{eq:tanSphere}).
\end{rem}

Since $n \ge 3$, the bottom line of~(\ref{eq:origintransform}) contains
at least $n$ equations.
We can assume $\bc_1,\ldots,\bc_n$ are linearly independent.
Then the matrix $M := (\bc_1, \ldots,\bc_n)^{\mathrm{T}}$ is invertible, and we can solve
the equations with indices $1, \ldots, n$ for~$p$:
 \begin{equation}\label{eq:pequation}
  p\ =\ \frac{1}{2\bv^2} M^{-1} 
    \left( \begin{array}{c}
      \bv^2 \bc_1^2 - [\bv \cdot \bc_1]^2 -\bv^2(r_1^2-r^2)\\
      \vdots \\
      \bv^2 \bc_n^2 - [\bv \cdot \bc_n]^2-\bv^2(r_n^2-r^2)\rule{0pt}{12.5pt}
    \end{array} \right).
 \end{equation}
Now substitute this expression for $p$ into the first and second equation
of the system~(\ref{eq:origintransform}), as well as into the equations 
for $i=n+1,\ldots,2n-3$, and then clear the denominators. 
This gives $n-1$ homogeneous equations in the coordinate $\bv$,
namely one cubic, one quartic, and $n-3$ quadrics.
By B\'{e}zout's Theorem, this means that if the system has only finitely
many solutions, then the number of solutions is bounded by
$3 \cdot 4 \cdot 2^{n-3} = 3 \cdot 2^{n-1}$, for $n \ge 3$. 
For small values of $n$, these values are shown in Table~$1$.
The value 12 for $n=3$ was computed in~\cite{MPT01},
and the values for $n=4,5,6$ were computed experimentally
in~\cite{sottile-macaulay-2001}.
\begin{table}[htb]
 \begin{center}
  \begin{tabular}{|c||c|c|c|c|c|}
  \hline
     $n$                 &  3 &  4 &  5 &  6  &   7\\ \hline
     maximum \# tangents & 12 & 24 & 48 & 96  & 192\\ 
   \hline
  \end{tabular}\smallskip
  \end{center}
 \label{ta:ndimvalues}
 \caption{Maximum number of tangents in small dimensions}
\end{table}

We simplify the cubic equation obtained by substituting~(\ref{eq:pequation})
into the equation $p\cdot v=0$ by expressing it in the basis 
$c_1, \ldots, c_n$.
Let the representation of $v$ in the basis $c_1, \ldots, c_n$
be 
\[
  v\ =\ \sum_{i=1}^n t_i c_i
\]
with homogeneous coordinates $t_1, \ldots, t_n$.
Further, let $c_1', \ldots, c_n'$ be a dual basis to $c_1, \ldots, c_n$;
i.e., let $c_1', \ldots, c_n'$ be defined by 
$c_i' \cdot c_j = \delta_{ij}$, where $\delta_{ij}$ denotes Kronecker's
delta function. 
By elementary linear algebra, we have $t_i = c_i' \cdot v$.

When expressing $p$ in this dual basis,
$p = \sum p_i' c_i'$, the third equation of~(\ref{eq:origintransform}) 
gives
\[
  p_i'\ =\ \frac{1}{v^2}
    \left(v^2 c_i^2 - [v \cdot c_i]^2 - v^2 (r_i^2 - r^2)\right) \, .
\]
Substituting this representation of $p$ into the equation
\[
  0\ =\ 2 v^2 (p \cdot v)
   \ =\ 2 v^2 \left(\sum_{i=1}^n p_i' c_i'\right) \cdot v 
   \ =\ 2 v^2 \sum_{i=1}^n p_i' t_i \, ,
\]
we obtain the cubic equation
\[
  \sum_{i=1}^n (v^2 c_i^2 - [v \cdot c_i]^2 - v^2 (r_i^2 - r^2)) t_i\ =\ 0 \, .
\]
In the case that all radii are equal, 
expressing $\bv^2$ in terms of the $t$-variables yields
\[
  \sum_{1 \le i \neq j \le n}
    \alpha_{ij} t_i^2 t_j + 
    \sum_{1 \le i < j < k \le n} 2 \beta_{ijk} t_i t_j t_k\ =\ 0 \, ,
\]
where
\begin{eqnarray*}
  \alpha_{ij} & = & (\text{vol}_2(\bc_i,\bc_j))^2\ =\ 
    \det \left( \begin{array}{cc}
      \bc_i\cdot \bc_i & \bc_i\cdot \bc_j  \\
       \bc_j\cdot \bc_i & \bc_j\cdot \bc_j 
  \end{array} \right), \\
  \beta_{ijk} & = &
    \det \left( \begin{array}{cc}
       \bc_i\cdot \bc_j  &  \bc_i\cdot \bc_k  \\
       \bc_k\cdot \bc_j  &  \bc_k\cdot \bc_k 
  \end{array} \right) +
    \det \left( \begin{array}{cc}
       \bc_i\cdot \bc_k  &  \bc_i\cdot \bc_j  \\
       \bc_j\cdot \bc_k  &  \bc_j\cdot \bc_j 
  \end{array} \right) \\
  & & +
    \det \left( \begin{array}{cc}
       \bc_j\cdot \bc_k  &  \bc_j\cdot \bc_i  \\
       \bc_i\cdot \bc_k  &  \bc_i\cdot \bc_i 
  \end{array} \right),
\end{eqnarray*}
and $\text{vol}_2(\bc_i,\bc_j)$ denotes the oriented area of the parallelogram
spanned by $\bc_i$ and $\bc_j$.
In particular, if $0c_1 \ldots c_n$ constitutes a regular simplex in 
$\mathbb{R}^n$,
then we obtain the following characterization.

\smallskip

\begin{thm} Let $n \ge 3$.
 If\/ $0\bc_1\ldots \bc_n$ is a regular simplex and all spheres have
 the same radius, 
 then the cubic equation expressed in the basis 
 $\bc_1, \ldots, \bc_n$ is equivalent to
 \begin{equation}
 \label{eq:cubicndimsimplex}
   \sum_{1 \le i \neq j \le n} t_i^2 t_j + 
   2 \sum_{1 \le i < j < k \le n} t_i t_j t_k\ =\ 0.
 \end{equation}
 For $n=3$, this cubic equation factors into three linear
terms;
 for $n \ge 4$ it is irreducible.
\end{thm}

\smallskip

\begin{proof} Let $e$ denote the edge length of the regular simplex.
Then the form of the cubic equation follows from computing
  $\alpha_{ij} = e^2 ( 1 \cdot 1 - 1/2 \cdot 1/2) = 3e^2/4$,
  $\beta_{ijk} = 3 e^2 (1/2 \cdot 1 - 1/2 \cdot 1/2) = 3e^2/4$.

Obviously, for $n=3$ the cubic polynomial
factors into $(t_1 + t_2)(t_1 + t_3)(t_2 + t_3)$ 
(cf.~\cite{schaal-85,MPT01}).
For $t \ge 4$, assume that there exists a factorization of the form
\[
  \left(t_1 + \sum_{i=2}^n \rho_i t_i\right)
  \left(\sum_{1 \le i \le j \le n} \sigma_{ij} t_i t_j\right)
\]
with $\sigma_{12} = 1$. 
Since~(\ref{eq:cubicndimsimplex}) does not
contain a monomial $t_i^3$, we have
either $\rho_i = 0$ or $\sigma_{ii} = 0$ for $1 \le i \le n$.

If there were more than one vanishing coefficient $\rho_i$, say
$\rho_i = \rho_j = 0$, then the monomials $t_i^2 t_j$ could not be
generated. So only two cases have to be investigated.

\smallskip

\noindent
\emph{Case~1}: $\rho_i \neq 0$ for $2 \le i \le n$. Then 
$\sigma_{ii} = 0$ for $1 \le i \le n$. Furthermore, $\sigma_{ij} = 1$
for $i \neq j$ and $\rho_i = 1$ for all $i$. Hence, the coefficient
of the monomial $t_1 t_2 t_3$ is 3, which 
contradicts~(\ref{eq:cubicndimsimplex}).

\smallskip

\noindent
\emph{Case~2}: There exists exactly one coefficient $\rho_i = 0$,
say, $\rho_4 = 0$. Then $\sigma_{11} = \sigma_{22} = \sigma_{33} = 0$,
$\sigma_{44} = 1$. Further, $\sigma_{ij} = 1$ for $1 \le i < j \le 3$  
and $\rho_i = 1$ for $1 \le i \le 3$. Hence, the coefficient of the monomial
$t_1 t_2 t_3$ is 3, which is again a contradiction.
\end{proof}
  
\section{Real Lines}

In the previous section, we have given the upper bound of
$3 \cdot 2^{n-1}$ for the number of complex solutions to the
tangent problem.
Now we complement this result by providing
a class of configurations leading to
$3 \cdot 2^{n-1}$ real common tangents. 
Hence, the upper bound is tight, and is achieved by real tangents.

There are no general techniques known to find and prove
configurations with a maximum number of real solutions in enumerative
geometry problems like the one studied here. 
For example, for the
classical enumerative geometry problem of 3264 conics tangent to five
given conics (dating back to Steiner in 1848~\cite{St1848} and solved by
Chasles in 1864~\cite{Ch1864}) the existence of five real conics with all 3264
\emph{real} was only recently established (\cite{rtv-97} and
\cite[{\S}7.2]{fulton-b96}). 

Our construction is based on the following geometric idea. For
4 spheres in $\mathbb{R}^3$ centered at the vertices $(1,1,1)^{\mathrm{T}}$,
$(1,-1,-1)^{\mathrm{T}}$, $(-1,1,-1)^{\mathrm{T}}$, $(-1,-1,1)^{\mathrm{T}}$ of a regular tetrahedron,
there are~\cite{MPT01} 
\begin{itemize}
  \item 3 different real tangents (of multiplicity~4) for radius 
        $r=\sqrt{2}$;
  \item 12 different real tangents for $\sqrt{2} < r < 3/2$;
  \item 6 different real tangents (of multiplicity~2) for $r=3/2$.
\end{itemize}
Furthermore, based on the explicit calculations in~\cite{MPT01}, 
it can be easily seen that the symmetry group of the
tetrahedron acts transitively on the tangents. By this symmetry
argument, all 12 tangents have the same 
distance $d$ from the origin. In order to construct a configuration
of spheres
with many common tangents, say, in $\mathbb{R}^4$, we embed the centers via
\[
  (x_1,x_2,x_3)^{\mathrm{T}}\ \longmapsto\ (x_1,x_2,x_3,0)^{\mathrm{T}}
\]
into $\mathbb{R}^4$ and place additional spheres with radius $r$ at
$(0,0,0,a)^{\mathrm{T}}$ and $(0,0,0,-a)^{\mathrm{T}}$ for some appropriate
value of $a$. 
If $a$ is chosen in such a way that the centers of the two additional 
spheres have distance $r$ from the above tangents,
then, intuitively, all common tangents to the six four-dimensional spheres
are located in the hyperplane $x_4 = 0$ and have multiplicity~2
(because of the two different
possibilities of signs when perturbing the situation).
By perturbing this configuration slightly, the tangents are no longer 
located in the hyperplane $x_4 = 0$, and therefore the double tangents are forced 
to split. The idea also generalizes to dimension $n \ge 5$.

Formally, suppose that the $2n-2$ spheres in $\mathbb{R}^n$ all have the same 
radius, $r$, and the first four have centers 
 \begin{eqnarray*}
 \bc_1 &:=& (\hh 1,\hh 1,\hh 1,\ 0,\ldots,0)^{\mathrm{T}}, \\
 \bc_2 &:=& (\hh 1,-1,-1,\ 0,\ldots,0)^{\mathrm{T}}, \\
 \bc_3 &:=& (-1,\hh 1,-1,\ 0,\ldots,0)^{\mathrm{T}}, \mbox{\quad and}\\
 \bc_4 &:=& (-1,-1,\hh 1,\ 0,\ldots,0)^{\mathrm{T}}
 \end{eqnarray*}
at the vertices of a regular tetrahedron inscribed in the 3-cube 
$(\pm1,\pm1,\pm1,0,\ldots,0)^{\mathrm{T}}$.
We place the subsequent centers at the points $\pm a\be_j$  for $j=4,5,\ldots,n$,
where $\be_1,\ldots,\be_n$ are the standard unit vectors in ${\mathbb R}^n$.

\begin{thm}\label{thm:real-cmplx}
Let $n \ge 4$, $r > 0$, $a > 0$, and $\gamma:=a^2(n-1)/(a^2+n-3)$.
If 
 \begin{equation}\label{eq:discrim}
  (r^2-3)\,(3-\gamma)\,(a^2-2)\,(r^2-\gamma)\,
  \left((3-\gamma)^2 +4\gamma - 4r^2\right)\ \neq\ 0\,,
 \end{equation}
%
%
then there are exactly 
$3\cdot 2^{n-1}$ different lines tangent to the $2n-2$ spheres. If 
 \begin{equation}\label{eq:ineqs}
   a^2 > 2,\quad  \gamma < 3,\quad \mbox{and}\quad
     \gamma\ <\ r^2\ <\
     \gamma + {\textstyle\frac{1}{4}}\left(3-\gamma\right)^2 \,,
 \end{equation}
then all these $3\cdot 2^{n-1}$ lines are real. 
Furthermore, this system of inequalities 
defines a nonempty subset of the $(a,r)$-plane. 
\end{thm}

Given values of $a$ and $r$ satisfying~(\ref{eq:ineqs}), we may scale the
centers and parameters by $1/r$ to obtain a configuration with unit spheres, 
proving Theorem~\ref{th:ndimnumber}~(b).

\begin{rem}
The set of values of $a$ and $r$ which give all solutions real is nonempty.
To show this, we calculate
 \begin{equation}\label{eq:gamma}
  \gamma\ =\ \frac{a^2(n-1)}{a^2+n-3}\ =\ 
    (n-1)\left(1-\frac{n-3}{a^2+n-3}\right)\,,
 \end{equation}
which implies that $\gamma$ is an increasing function of $a^2$.
Similarly, set $\delta:=\gamma+(3-\gamma)^2/4$, the upper bound
for $r^2$.
Then
$$
  \frac{d}{d\gamma}\;\delta\ =\ 
    \frac{d}{d\gamma}\left(\frac{\gamma + (3-\gamma)^2}{4}\right)p\ =\ 
    1 + \frac{\gamma-3}{2}\,,
$$
and so $\delta$ is an increasing function of $\gamma$ when $\gamma>1$.
When $a^2=2$, we have $\gamma=2$; so 
$\delta$ is an increasing function of $a$ in the region $a^2>2$.
Since when $a=\sqrt{2}$, 
we have $\delta=\frac{9}{4}>\gamma$, the region defined
by~(\ref{eq:ineqs}) is nonempty.

Moreover, we remark that the region is qualitatively different in the 
cases $n=4$ and $n \ge 5$. For $n=4$, $\gamma$ satisfies $\gamma<3$ 
for any $a>\sqrt{2}$. Hence, $\delta<3$ and $r<\sqrt{3}$.
Thus the maximum value of 24 real lines may be obtained for arbitrarily
large $a$. In particular, we may choose the two spheres with 
centers $\pm ae_4$ disjoint from the first four spheres. 
Note, however, that the first four spheres do meet, since we have
$\sqrt{2}<r<\sqrt{3}$.

For $n\geq 5$, there is an upper bound to $a$.
The upper and lower bounds for $r^2$ coincide when $\gamma=3$;
so we always have $r^2<3$.
Solving $\gamma=3$ for $a^2$, we obtain $a^2<3(n-3)/(n-4)$.
When $n=5$, Figure~\ref{fig:realShade} displays the discriminant locus
(defined by~(\ref{eq:discrim})) and
shades the region consisting of values of $a$ and $r$ for which all solutions
are real. 

\ifpictures
\begin{figure}[htb]
$$
  \setlength{\unitlength}{1.32pt}  
  \begin{picture}(308,145)(3,0)
    \put(0,0){\epsfxsize=403pt\epsffile{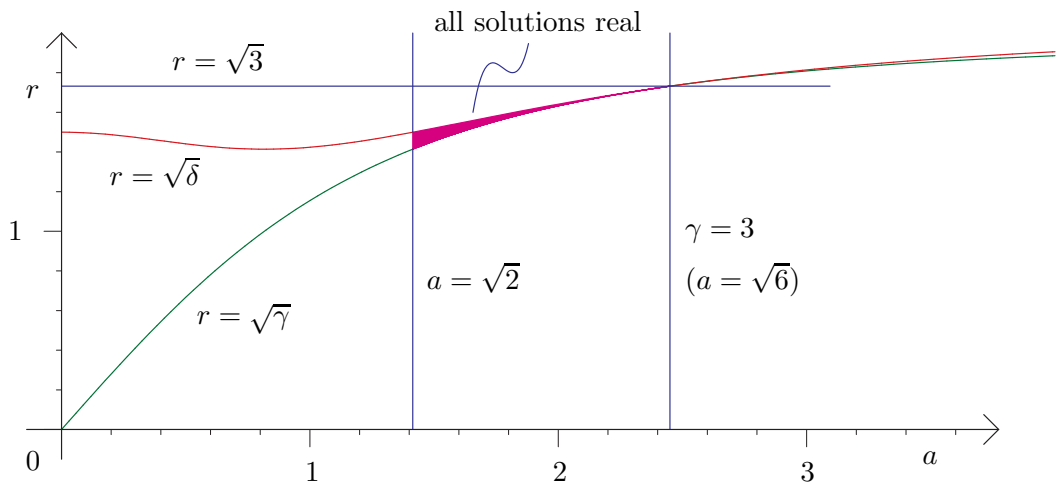}}  
    \put(265,10){$a$}
    \put(8,115){$r$}

    \put(8, 8){$0$}
    \put(3,74){$1$}

    \put( 88,5){$1$}
    \put(159,5){$2$}
    \put(230,5){$3$}

    \put(50,123){$r=\sqrt{3}$}
    \put(32, 90){$r=\sqrt{\delta}$}
    \put(57, 50){$r=\sqrt{\gamma}$}

    \put(123,60){$a=\sqrt{2}$}
    \put(197,75){$\gamma=3$}
    \put(197,60){$(a=\sqrt{6})$}

    \put(125,134){all solutions real}

  \end{picture}
$$
\caption{Discriminant locus and values of $a,r$ giving all solutions
real\label{fig:realShade}} 
\end{figure}
\fi

\end{rem}

\noindent
\emph{Proof of Theorem~\ref{thm:real-cmplx}}.
We prove Theorem~\ref{thm:real-cmplx} by treating $a$ and $r$ as parameters
and explicitly solving the resulting system 
of polynomials in the coordinates 
$(\bp,\bv)\in{\mathbb C}^n\times{\mathbb P}_{\mathbb C}^{n-1}$ 
for lines in ${\mathbb C}^n$.
This shows that there are $3\cdot 2^{n-1}$ {\it complex} lines tangent to the
given spheres, for the values of the parameters $(a,r)$ given in 
Theorem~\ref{thm:real-cmplx}.
The inequalities~(\ref{eq:ineqs}) describe the parameters for which all
solutions are real. 
\medskip

First consider the equations~(\ref{eq:tanSphere}) for the line to be tangent
to the spheres with centers $\pm a \be_j$ and radius $r$:
 \begin{eqnarray*}
  \bv^2\bp^2-2 a \bv^2p_j + a^2\bv^2 - a^2v_j^2-r^2\bv^2&=&0,\\
  \bv^2\bp^2+2 a \bv^2p_j + a^2\bv^2 - a^2v_j^2-r^2\bv^2&=&0.
 \end{eqnarray*}
Taking their sum and difference (and using $a\bv^2\neq0$), we obtain
 \begin{eqnarray}
   p_j&=&0, \hspace{2.64cm}\qquad 4 \le j \le n, \label{eq:pj=0}\\
   a^2v_j^2  &=&(\bp^2+a^2-r^2)\bv^2, \qquad\ 4 \le j \le n. \label{eq:vj}
 \end{eqnarray}
Subtracting the equations~(\ref{eq:tanSphere}) for the centers $c_1, \ldots, c_4$
pairwise gives
\[
  4 v^2 (p_2 + p_3)\ =\ -4 (v_1 v_3 + v_1 v_2)
\]
(for indices 1,2) and analogous equations. Hence,
\[
  p_1\ =\ - \frac{v_2 v_3}{v^2}, \qquad
  p_2\ =\ - \frac{v_1 v_3}{v^2}, \qquad
  p_3\ =\ - \frac{v_1 v_2}{v^2}.
\]
Further, $p \cdot v = 0$ implies $v_1 v_2 v_3 = 0$. 
Thus we have 3 symmetric cases.
We treat one, assuming that $v_1 = 0$.
Then we obtain
\[
  p_1\ =\ - \frac{v_2 v_3}{v^2}, \qquad p_2 = p_3 = 0.
\]
Hence, the tangent equation~(\ref{eq:tanSphere}) for the first sphere becomes
$$
  \bv^2p_1^2-2\bv^2p_1 + 3\bv^2-(v_2+v_3)^2-r^2\bv^2\ =\ 0\,.
$$
Using $0 = \bv^2p_1+v_2v_3$, we obtain
 \begin{equation}\label{eq:redE1}
   v_2^2+v_3^2\ =\ \bv^2(p_1^2+3 -r^2)\,.
 \end{equation}
The case $j=4$ of (\ref{eq:vj}) gives $a^2v_4^2=\bv^2(p_1^2+a^2-r^2)$,
since $p_2=p_3=0$.
Combining these, we obtain
$$
  v_2^2+v_3^2\ =\ a^2 v_4^2 + \bv^2(3-a^2)\,.
$$
Using $\bv^2=v_2^2+v_3^2+(n-3)v_4^2$ yields
\[
  (a^2-2)(v_2^2+v_3^2)\ =\ v_4^2(3(a^2+n-3)-a^2(n-1)).
\]
We obtain
 \begin{equation}\label{eq:v234}
    (a^2-2)(v_2^2+v_3^2)\ =\ v_4^2(a^2+n-3)(3-\gamma)\,,
 \end{equation}
where $\gamma=a^2(n-1)/(a^2+n-3$).

Note that $a^2+n-3>0$ since $n>3$.
If neither $3-\gamma$ nor $a^2-2$ are zero, then we may use this to compute
 \begin{eqnarray*}
  (a^2+n-3)(3-\gamma)\bv^2 &=& [(a^2+n-3)(3-\gamma)+(n-3)(a^2-2)](v_2^2+v_3^2)\\
                      &=& (a^2+n-3)(v_2^2+v_3^2)\, ,
 \end{eqnarray*}
and so 
 \begin{equation}\label{eq:RealOne}
  (3-\gamma)\bv^2\ =\ v_2^2+v_3^2\,.
 \end{equation}
Substituting~(\ref{eq:RealOne}) into~(\ref{eq:redE1}) and dividing by $\bv^2$
gives 
 \begin{equation}\label{eq:RealTwo}
  p_1^2 \ =\ r^2-\gamma\,.
 \end{equation}
Combining this with $\bv^2p_1+v_2v_3=0$, we obtain
 \begin{equation}\label{eq:RealThree}
  p_1(v_2^2+v_3^2) + (3-\gamma)v_2v_3\ =\ 0\,.
 \end{equation}\medskip
Summarizing, we have $n$ linear equations
$$
  v_1\ =\ p_2\ =\ p_3\ =\ p_4\ =\ \cdots\ =\ p_n\ =\ 0\,,
$$
and $n-4$ simple quadratic equations 
$$
  v_4^2\ =\ v_5^2\ =\ \cdots\ =\ v_n^2\,,
$$
and the three more complicated quadratic
equations,~(\ref{eq:v234}),~(\ref{eq:RealTwo}), and~(\ref{eq:RealThree}).
\smallskip

We now solve these last three equations.
We solve~(\ref{eq:RealTwo}) for $p_1$, obtaining
$$
  p_1\ =\ \pm\sqrt{r^2-\gamma}\,.
$$
Then we solve~(\ref{eq:RealThree}) for $v_2$ and use~(\ref{eq:RealTwo}), obtaining
$$
   v_2\ =\ -\frac{3-\gamma\pm\sqrt{(3-\gamma)^2-4(r^2-\gamma)}}{2 p_1}\, v_3\,.
$$
Finally,~(\ref{eq:v234}) gives
$$
   v_4\sqrt{a^2+n-3}\ =\ \pm\sqrt{\frac{a^2-2}{3-\gamma}(v_2^2+v_3^2)}\,.
$$
Since $v_3 = 0$ would imply $v=0$ and hence contradict 
$v \in {\mathbb P}^{n-1}_{\mathbb C}$, we see that $v_3 \neq 0$. Thus we can
conclude that when none of the following expressions
$$
  r^2-3\,,\;\ 3-\gamma\,,\;\ a^2-2\,,\;\ r^2-\gamma\,,\; \ 
  (3-\gamma)^2 +4\gamma - 4r^2
$$
vanish, there are $8=2^3$ different solutions to the last 3 equations.
For each of these, the simple quadratic equations give $2^{n-4}$ solutions;
so we see that the case $v_1 = 0$ contributes $2^{n-1}$ different solutions,
each of them satisfying $v_2 \neq 0$, $v_3 \neq 0$.
Since there are three symmetric cases, we obtain $3\cdot2^{n-1}$ solutions in
all, as claimed.
\bigskip

We complete the proof of Theorem~\ref{thm:real-cmplx} and determine
which values of the parameters $a$ and $r$ give all these lines real.
We see that
\begin{enumerate}
  \item[(1)] $p_1$ is real if $r^2-\gamma >0$.
  \item[(2)] Given that $p_1$ is real, $v_2/v_3$ is real if 
        $(3-\gamma)^2+4\gamma-4r^2>0$.
  \item[(3)] Given this, $v_4/v_3$ is real if $(a^2-2)/(3-\gamma)>0$.
\end{enumerate}

Suppose the three inequalities above are satisfied.
Then all solutions are real, and~(\ref{eq:RealOne}) implies that $3-\gamma>0$,
and so we also have $a^2-2>0$.
This completes the proof of Theorem~\ref{thm:real-cmplx}.  \qed

\section{Affinely Dependent Centers}
In our derivation of the B\'ezout number $3\cdot 2^{n-1}$ of common tangents 
for Theorem~\ref{th:ndimnumber}, it was crucial that the centers of the
spheres affinely spanned ${\mathbb R}^n$.
Also, the construction in Section 3 of configurations with $3\cdot 2^{n-1}$ 
real common tangents had centers affinely spanning ${\mathbb R}^n$.
When the centers are affinely dependent, we prove the following result.

\begin{thm}\label{thm:af-dep} For $n \ge 4$,
 there are $3\cdot 2^{n-1}$ complex common tangent lines to $2n-2$ spheres
 whose centers are affinely dependent, but otherwise general.
 There is a choice of such spheres with $2^n$ real common tangent lines.
\end{thm}

\begin{rem}
 Theorem~\ref{thm:af-dep} extends the results of~\cite[Section 4]{MPT01},
 where it is shown that when $n=3$, there are 12 complex common tangents.
 Megyesi~\cite{Me01} has shown that there is a configuration with
 12 real common tangents, but that the number of tangents is bounded
 by 8 for the case of unit spheres.
 For $n \ge 4$, we are unable either to find a configuration of spheres
 with affinely dependent centers and equal radii having more than $2^n$ real
 common tangents, 
 or to show that the maximum number of real common tangents is less than
 $3\cdot 2^{n-1}$.
 Similar to the case $n=3$, it might be possible that the case of unit 
 spheres and the case
 of spheres with general radii might give different maximum numbers.
\end{rem}

\begin{rem}
 Megyesi~\cite{Me02} recently showed that there are $2n-2$ spheres with
 affinely dependent centers having all $3\cdot 2^{n-1}$ common tangents real.
 Furthermore,  all but one of the spheres in his construction have equal radii.
\end{rem}

 By Theorem~\ref{th:ndimnumber}, $3\cdot 2^{n-1}$ is the upper bound for
the number of complex common tangents to spheres with affinely dependent
centers. 
Indeed, if there were a configuration with more common tangents, then---since
the system is a complete intersection---perturbing the centers would give a
configuration with affinely independent centers and more common tangent
lines than allowed by Theorem~\ref{th:ndimnumber}.

 By this discussion, to prove Theorem~\ref{thm:af-dep} it suffices to give
$2n-2$ spheres with affinely dependent centers having $3\cdot 2^{n-1}$ complex
common tangents and also such a configuration of $2n-2$ spheres with $2^n$
real common tangents. 
For this, we use spheres with equal radii whose centers are the vertices
of a perturbed crosspolytope in a hyperplane.
We work with the notation of Sections 2 and 3.

Let $a\neq -1$ and suppose we have spheres with equal radii $r$ and centers at
the points
$$
  ae_2,\ \,-e_2,\quad\mbox{and}\quad \pm e_j, \quad \mbox{for}\ 3\leq j\leq n\,.
$$
Then we have the equations
 \begin{eqnarray}
  p\cdot v&=&0,\label{eq:first}\\
  f\ :=\ v^2(p^2-2ap_2+a^2-r^2)-a^2v_2^2&=&0,\\
  g\ :=\hspace{2.24em} v^2(p^2+2p_2+1-r^2)-v_2^2&=&0,\label{eq:XX}\\
     v^2(p^2\pm2p_j+1-r^2)-v_j^2&=&0, \qquad 3\leq j\leq n\,.\label{eq:j}
 \end{eqnarray}

As in Section 3, the sum and difference of the equations~(\ref{eq:j})
for the spheres with centers $\pm e_j$ give
$$
 \begin{array}{rcl}
  p_j&=&0,\\
  v^2(p^2+1-r^2)&=&v_j^2.\rule{0pt}{14pt}
 \end{array}\qquad 3\leq j\leq n\,.
$$
Thus we have the equations
 \begin{equation}\label{eq:jj}
  \begin{array}{c}
    p_3\ =\ p_4\ =\ \cdots\ =\ p_n\ =\ 0,\\
    v_3^2\ =\ v_4^2\ =\ \cdots\ =\ v_n^2. \rule{0pt}{15pt}
  \end{array}
 \end{equation}
Similarly, we have
 \begin{eqnarray*}
    f+ag&=& (1+a)\left[v^2(p^2-r^2+a) - av_2^2\right] = 0,\\
    f-a^2g&=& (1+a)v^2\left[ (1-a)(p^2-r^2) + 2ap_2\right] = 0.
 \end{eqnarray*}
As before, $v^2\neq 0$:
If $v^2=0$, then~(\ref{eq:XX}) and~(\ref{eq:j}) imply that $v_2=\cdots=v_n=0$.
With $v^2=0$, this implies that $v_1=0$ and hence $v=0$,
contradicting $\bv\in{\mathbb P}^{n-1}_{\mathbb C}$.
By~(\ref{eq:jj}), we have $p^2=p_1^2+p_2^2$, and so we obtain the system of
equations in the variables $p_1,p_2,v_1,v_2,v_3$:
 \begin{equation}\label{eq:sm-sys}
  \begin{array}{r}
    p_1v_1+p_2v_2\ =\ 0,\\
    (1-a)(p_1^2+p_2^2-r^2) + 2ap_2\ =\ 0\rule{0pt}{15pt},\\
    v^2(p_1^2+p_2^2-r^2+a) - av_2^2\ =\ 0\rule{0pt}{15pt},\\
    v^2(p_1^2+p_2^2-r^2+1) - v_3^2\ =\ 0\rule{0pt}{15pt}.
  \end{array}
 \end{equation}
(For notational sanity, we do not yet make the substitution
$v^2=v_1^2+v_2^2+(n-2)v_3^2$.) 

We assume that $a\neq1$ and will treat the case $a=1$ at the end of this
section. 
Using the second equation of~(\ref{eq:sm-sys}) to cancel the terms 
$v^2(p_1^2+p_2^2)$ from the third equation and dividing the result by $a$, we
can solve for $p_2$:
 $$
  p_2\ =\ \frac{(1-a)(v^2-v_2^2)}{2v^2}\,.
 $$
If we substitute this into the first equation of~(\ref{eq:sm-sys}), we may
solve for $p_1$:
 $$
  p_1\ =\ -\frac{(1-a)(v^2-v_2^2)v_2}{2v^2v_1}\,.
 $$

Substitute these into the second equation of~(\ref{eq:sm-sys}), clear
the denominator $(4v_1^2v^4)$, and remove the common factor $(1-a)$ to obtain
the sextic 
 \begin{equation}\label{eq:sextic}
  (1-a)^2(v_1^2+v_2^2)(v^2-v_2^2)^2 \,-\, 4r^2v_1^2v^4 \,+\, 
     4av_1^2v^2(v^2-v_2^2)\ =\ 0\,.
 \end{equation}
Subtracting the third equation of~(\ref{eq:sm-sys}) from the fourth equation
and recalling that $v^2=v_1^2+v_2^2+(n-2)v_3^2$, we obtain
the quadratic equation
 \begin{equation}\label{eq:quadric}
  (1-a)v_1^2+v_2^2+ \left[(n-3)-a(n-2)\right]v_3^2\ =\ 0\,.
 \end{equation}
Consider the system consisting of the two equations~(\ref{eq:sextic})
and~(\ref{eq:quadric}) in the homogeneous coordinates $v_1,v_2,v_3$.
Any solution to this system gives a solution to the
system~(\ref{eq:sm-sys}), and thus gives $2^{n-3}$ solutions to the original
system~(\ref{eq:first})--(\ref{eq:j}).

These last two equations~(\ref{eq:sextic}) and~(\ref{eq:quadric}) are
polynomials in the squares of the variables $v_1^2, v_2^2, v_3^2$. 
If we substitute $\alpha=v_1^2, \beta=v_2^2$, and $\gamma=v_3^2$, then we have
a cubic and a linear equation, and any solution $\alpha,\beta,\gamma$ to these
with nonvanishing coordinates gives 4 solutions to the
system~(\ref{eq:sextic}) and~(\ref{eq:quadric}):
$(v_1,v_2,v_3)^{\mathrm{T}}:=(\alpha^{1/2},\pm\beta^{1/2},\pm\gamma^{1/2})^{\mathrm{T}}$, as
$v_1,v_2,v_3$ are homogeneous coordinates. 

Solving the linear equation in $\alpha,\beta,\gamma$ for $\beta$ and
substituting into the cubic equation gives a homogeneous cubic in 
$\alpha$ and $\gamma$ whose coefficients are polynomials in
$a,n,r$
\silentfootnote{Maple V.5 code verifying this and other
explicit calculations presented in this manuscript is available at 
{\tt www.math.umass.edu/\~{}sottile/pages/spheres}.}. 
The discriminant of this cubic is a polynomial with 
integral coefficients of degree 16 in the variables $a,n,r$ having 116 terms.
Using a computer algebra system, it can be verified that
this discriminant is irreducible over the rational numbers.
Thus, for any fixed integer $n\geq 3$, the discriminant is a nonzero
polynomial in $a,r$.
This implies that the cubic has 3 solutions for general $a,r$ and any integer
$n$.
Since the coefficients of this cubic similarly are nonzero polynomials for any
$n$, the solutions $\alpha,\beta,\gamma$ will be nonzero for general $a,r$
and any $n$.
We conclude:
$$
  \mbox{\begin{minipage}[c]{4.3in} 
    For any integer $n\geq 3$ and general $a,r$, there will be
    $3\cdot 2^{n-1}$ complex common tangents to spheres of radius $r$ with
    centers
    $$
       ae_2,\ \,-e_2,\quad\mbox{and}\quad \pm e_j, \
        \quad \mbox{for}\ 3\leq j\leq n\,.
    $$
  \end{minipage}}
$$

We return to the case when $a=1$, i.e., the centers are the
vertices of the crosspolytope $\pm e_j$ for $j=2,\ldots,n$.
Then our equations~(\ref{eq:jj}) and~(\ref{eq:sm-sys}) become
 \begin{equation}\label{eq:2n}
  \begin{array}{r}
    p_2\ =\ p_3\ =\ \cdots\ =\ p_n\ =\ 0,\\
    v_2^2\ =\ v_3^2\ =\ \cdots\ =\ v_n^2, \rule{0pt}{15pt}\\
    p_1v_1\ =\ 0,\rule{0pt}{15pt}\\
    v^2(p_1^2-r^2+1) - v_2^2\ =\ 0.\rule{0pt}{15pt}
  \end{array}
 \end{equation}
As before, $v^2=v_1^2+(n-1)v_2^2$.
We solve the last two equations.
Any solution they have (in ${\mathbb C}^1\times{\mathbb P}^1_{\mathbb C}$) gives rise to
$2^{n-2}$ solutions, by the second list of equations
$v_3^2=\cdots=v_n^2$.
By the penultimate equation $p_1v_1=0$, one of $p_1$ or $v_1$ vanishes.
If $v_1=0$, then the last equation becomes
$$
  (n-1)v_2^2(p_1^2-r^2+1)\ =\ v_2^2\,.
$$
Since $v_2=0$ implies $v^2=0$, we have $v_2\neq 0$ and so we may divide by
$v_2^2$ and solve for $p_1$ to obtain
$$
  p_1\ =\ \pm\sqrt{r^2-1+\frac{1}{n-1}}\,.
$$
If instead $p_1=0$, then we solve the last equation to obtain
$$
  \frac{v_1}{v_2}\ =\ \pm\sqrt{\frac{1}{1-r^2}+1-n}\,.
$$

Thus for general $r$, there will be $2^n$ common tangents to the spheres with
radius $r$ and centers $\pm e_j$ for $j=2,\ldots,n$.
We investigate when these are real.

We will have $p_1$ real when $r^2\ >\ 1-1/(n-1)$.
Similarly, $v_1/v_2$ will be real when  $1/(1-r^2)\ >\ n-1$.
In particular, $1-r^2>0$ and so $1>r^2$. Using this we get
$$
  1-r^2\ <\ \frac{1}{n-1}\qquad\mbox{so that}\qquad
   r^2\ >\ 1-\frac{1}{n-1}\,,
$$
which we previously obtained.

We conclude that there will be $2^n$ real common tangents to the spheres with
centers  $\pm e_j$ for $j=2,\ldots,n$ and radius $r$ when
$$
   \sqrt{1-\frac{1}{n-1}}\ <\ r\ <\ 1 \,.
$$
This concludes the proof of Theorem~\ref{thm:af-dep}.


\section{Lines Tangent to Quadrics}
Suppose that in our original question we ask for common tangents to
ellipsoids, or to more general quadric hypersurfaces.
Since all smooth quadric hypersurfaces are projectively equivalent, a natural
setting for this question is the following:

``How many common tangents are there to $2n-2$ general quadric hypersurfaces
in (complex) projective space ${\mathbb P}^n_{\mathbb C}$?''

\begin{thm}\label{thm:bezout}
 There are at most
$$
  2^{2n-2}\cdot\frac{1}{n}\binom{2n-2}{n-1}
$$
isolated common tangent lines to $2n-2$ quadric hypersurfaces 
in ${\mathbb P}^n_{\mathbb C}$.
\end{thm}

\begin{proof}
 The space of lines in ${\mathbb P}^n_{\mathbb C}$ 
is the Grassmannian of 2-planes in
 ${\mathbb C}^{n+1}$.
 The Pl\"ucker embedding~\cite{MR48:2152} realizes this as a projective
 subvariety of ${\mathbb P}_{\mathbb C}^{\binom{n+1}{2}-1}$ of degree
$$
  \frac{1}{n}\binom{2n-2}{n-1}\,.
$$
The theorem follows from the refined B\'ezout 
theorem~\cite[\S 12.3]{Fu84a} and from the fact that the condition
for a line to be tangent to a quadric hypersurface is a homogeneous quadratic
equation in the Pl\"ucker coordinates for 
lines~\cite[\S 5.4]{sottile-macaulay-2001}.
\end{proof}

In Table~\ref{ta:quadvalues}, we compare the upper
bound of Theorem~\ref{thm:bezout} for the number of
lines tangent to $2n-2$ quadrics to the number of lines tangent to $2n-2$
spheres of Theorem~\ref{th:ndimnumber}, for small values of $n$.
\begin{table}[htb]
 \begin{center}
  \begin{tabular}{|c||c|c|c|c|c|}
  \hline
     $n$              &  3 &  4 &  5 &  6  &   7\\ \hline
     \# for spheres  & 12 & 24 & 48 & 96  & 192\\ \hline
     \# for quadrics &32 &320 &3580 &43008&540672\\ 
   \hline
  \end{tabular}\smallskip
  \end{center}
 \label{ta:quadvalues}
 \caption{Maximum number of tangents in small dimensions}
\end{table}
The bound of 32 tangent lines to 4 quadrics in ${\mathbb P}^3_{\mathbb C}$
is sharp, even under the restriction to real quadrics and real tangents
\cite{sottile-theobald-progress}.
In a computer calculation, we found 320 lines in ${\mathbb P}^4_{\mathbb C}$ 
tangent to 6 general quadrics; thus, the upper bound of
Theorem~\ref{thm:bezout} is sharp also for $n=4$, and indicating that it is
likely sharp for $n > 4$.
The question arises: what is the source of the huge discrepancy between the
second and third rows of Table~\ref{ta:quadvalues}?

Consider a sphere in affine $n$-space
$$
  (x_1-c_1)^2 + (x_2-c_2)^2 + \cdots + (x_n-c_n)^2\ =\ r^2\,.
$$
Homogenizing this with respect to the new variable $x_0$, we obtain
$$
  (x_1-c_1x_0)^2 + (x_2-c_2x_0)^2 + \cdots + (x_n-c_nx_0)^2\ 
      =\ r^2x_0^2\,.
$$
If we restrict this sphere to the hyperplane at infinity, setting $x_0=0$, we
obtain
 \begin{equation}\label{eq:imaginary}
  x_1^2+x_2^2+\cdots+x_n^2\ =\ 0\,,
 \end{equation}
the equation for an imaginary quadric at infinity.
We invite the reader to check that every line at infinity tangent to this
quadric is tangent to the original sphere.

Thus the equations for lines in ${\mathbb P}^n_{\mathbb C}$ tangent to $2n-2$ spheres
define the $3\cdot 2^{n-1}$ lines we computed in Theorem~\ref{th:ndimnumber},
as well as this excess component of lines at infinity tangent to the imaginary 
quadric~(\ref{eq:imaginary}).
Thus, this excess component contributes some portion of the
B\'ezout number of Theorem~\ref{thm:bezout} to the total number of lines.
Indeed, when $n=3$, Aluffi and Fulton~\cite{aluffi-fulton-2001} have given a
careful argument that this excess component contributes 20, which implies
that there are $32-20=12$ isolated common tangent lines to 4 spheres in 3-space,
recovering the result of~\cite{MPT01}.

The geometry of that calculation is quite interesting.
Given a system of equations on a space (say the Grassmannian) whose set of
zeroes has a positive-dimensional excess component, one method to compute the
number of isolated solutions is to first modify the underlying space by
blowing up the 
excess component and then compute the number of solutions on this new space.
In many cases, the equations on this new space have only isolated solutions.
However, for this problem of lines tangent to spheres, the equations on the
blown up space will \emph{still} have an excess intersection and a further
blow-up is required. 
This problem of lines tangent to 4 spheres in projective 3-space is by far the
simplest enumerative geometric problem with an excess component of zeroes
which requires two blow-ups (technically speaking, blow-ups along smooth
centers) to resolve the excess zeroes.

It would be interesting to understand the geometry also when 
$n > 3$.
For example, how many blow-ups are needed to resolve the excess component?
\bigskip

Since all smooth quadrics are projectively equivalent,
Theorem~\ref{th:ndimnumber} has the following implication for this problem of
common tangents to projective quadrics.

\begin{thm}\label{thm:UB}
 Given $2n-2$ quadrics in ${\mathbb P}^n_{\mathbb C}$ 
whose intersection with a fixed
 hyperplane is a given smooth quadric $Q$, but are otherwise general, there will
 be at most $3\cdot 2^{n-1}$ isolated lines in ${\mathbb P}^n_{\mathbb C}$ 
tangent to each
 quadric. 
\end{thm}

When the quadrics are all real, we ask: how many of these  $3\cdot 2^{n-1}$
common isolated tangents can be real?
This question is only partially answered by Theorem~\ref{th:ndimnumber}.
The point is that projective real quadrics are classified up to real
projective transformations by the absolute value of the signature of the
quadratic forms on ${\mathbb R}^{n+1}$ defining them.
Theorem~\ref{th:ndimnumber} implies that all lines can be real when the shared
quadric $Q$ has no real points (signature is $\pm n$).
In~\cite{sottile-macaulay-2001}, it is shown that when $n=3$, each of the five
additional cases concerning nonempty quadrics can have all 12 lines real.
\smallskip

Recently, Megyesi~\cite{Me02} has largely answered this question.
Specifically, he showed that, for any nonzero real numbers 
$\lambda_3,\ldots,\lambda_n$, there are $2n-2$ quadrics of the form
$$
  (x_1-c_1)^2 + (x_2-c_2)^2 + \sum_{j=3}^n \lambda_j(x_j-c_j)^2
  \ =\ R
$$
having all $3\cdot 2^{n-1}$ tangents real.
These all share the same quadric at infinity
$$
  x_1^2 + x_2^2 + \lambda_3 x_3^2 + \cdots + \lambda_n x_n^2\ =\ 0\,,
$$
and thus the upper bound of Theorem~\ref{thm:UB} is attained, when the shared
quadric is this quadric.
\bigskip

\noindent
{\bf Acknowledgments:} The authors would like to thank
I.~G.~Macdonald for pointing out a
simplification in Section~2, as well as Gabor Megyesi and an
unkwown referee for their useful suggestions.

\providecommand{\bysame}{\leavevmode\hbox to3em{\hrulefill}\thinspace}
\providecommand{\MR}{\relax\ifhmode\unskip\space\fi MR }
\providecommand{\MRhref}[2]{%
  \href{http://www.ams.org/mathscinet-getitem?mr=#1}{#2}
}
\providecommand{\href}[2]{#2}

\end{document}

\bibliographystyle{amsplain}
\bibliography{bibl}

\end{document}